\documentclass[a4paper,12pt]{article}
\usepackage{latexsym,amsmath,amsthm,amssymb}
\usepackage{a4wide}
\usepackage{hyperref}
\usepackage{marginnote}
\usepackage{color}
\hypersetup{
pdftitle={Improved $L_p-$mixed volume inequality for convex bodies},    
pdfauthor={Van Hoang Nguyen},
colorlinks = true,
linkcolor = magenta,
citecolor = blue,
}

\theoremstyle{plain}
\newtheorem{theorem}{Theorem}[section]

\theoremstyle{definition}

\theoremstyle{remark}

\renewcommand{\thefootnote}{\arabic{footnote}}

\def\R{\mathbb R}

\def\Om{\Omega}
\def\be{\beta}
\def\ga{\gamma}
\def\de{\delta}
\def\De{\Delta} 

\def\si{\sigma}
\def\lam{\lambda}
\def\ep{\epsilon}
\def\la{\langle} 
\def\ra{\rangle} 
\def\lt{\left}
\def\rt{\right}

\def\mK{\mathcal K}

\def\i0i{\int_0^\infty}

\numberwithin{equation}{section}

\title{Improved $L_p-$mixed volume inequality for convex bodies}
\author{Van Hoang Nguyen\footnote{
School of Mathematical Sciences, Tel Aviv University, Tel Aviv 69978, \textsc{Israel}.}}

\begin{document}
\maketitle


\renewcommand{\thefootnote}{}

\footnote{Email: vanhoang0610@yahoo.com}

\footnote{Supported by a grant from the European Research Council (grant number $305629$)}

\footnote{2010 \emph{Mathematics Subject Classification\text}: 26D15, 52A20, 52A39, 52A40.}

\footnote{\emph{Key words and phrases\text}: $L_p-$mixed volume, Brunn-Minkowski-Firey theory, stability, Jensen inequality.}

\renewcommand{\thefootnote}{\arabic{footnote}}
\setcounter{footnote}{0}

\begin{abstract} 
A sharp quantitative version of the $L_p-$mixed volume inequality is established. This is achieved by exploiting an improved Jensen inequality. This inequality is a generalization of Pinsker-Csisz\'ar-Kullback inequality for the Tsallis entropy. Finally, a sharp quantitative version of the $L_p-$Brunn-Minkowski inequality is also proved as a corollary.
\end{abstract}

\section{Introduction}
Throughout this paper, a convex body is a compact convex subset with nonempty interior in $\R^n$. It is well known that a convex body $K$ is uniquely determined by its support function defined by 
$$h(K,x) = \max\{\la x,y\ra\,:\, y\in K\},$$
where $\la\cdot,\cdot\ra$ is the Euclidean scalar product in $\R^n$.

Let $\mK_0^n$ denote the set of convex bodies containing the origin in its interior. For each $p\geq 1$, convex bodies $K,L\in \mK_0^n$ and $\ep > 0$, the Minkowski-Firey $L_p-$combination $K+_p\ep\cdot L$ which was firstly introduced and studied by Firey \cite{Firey62} is the convex body whose support function is given by 
$$h(K+_p\ep\cdot L, x) = \lt(h(K,x)^p + \ep h(L,x)^p\rt)^{\frac 1p},\quad x\in \R^n.$$ 
For $K,L\in \mK_0^n$, their $L_p-$mixed volume $V_p(K,L)$ is defined by
$$V_p(K,L) = \frac pn \lim\limits_{\ep\to 0} \frac{V(K+_p\ep\cdot L) -V(K)}\ep,$$
where $V(K)$ stands for the volume of $K$. It was shown by Lutwak (see \cite{Lut91}) that this limit exists and has the following integral representation
\begin{equation}\label{eq:intrepre}
V_p(K,L) = \frac1n \int_{S^{n-1}} \lt(\frac{h(L,u)}{h(K,u)}\rt)^p h(K,u) dS(K,u),
\end{equation}
where $S(K,\cdot)$ denotes surface area measure of $K$. It is obvious that $V_p(K,K) =V(K)$ for any $K\in \mK_0^n$.

The $L_p-$mixed volume satisfies the following inequality (see \cite[Theorem $1.2$]{Lut91})
\begin{equation}\label{eq:Lpmixed}
V_p(K,L) \geq V(K)^{1-\frac pn} V(L)^{\frac pn}, \quad \forall\, K,L \in \mK_0^n,
\end{equation}
with equality holds if and only if $K,L$ are dilates. It is also proved in \cite{Lut91} that 
\begin{equation}\label{eq:extensionBM}
V(K+_p L)^{\frac pn} \geq V(K)^{\frac pn} + V(L)^{\frac pn},\quad\forall\, K,L \in \mK_0^n,
\end{equation}
with the same equality condition as above. The inequality \eqref{eq:extensionBM} is Firey's extension of the famous Brunn-Minkowski inequality corresponding to $p=1$.

The motivation of this paper is to study the quantitative improvement version of the inequality \eqref{eq:Lpmixed} and \eqref{eq:extensionBM}. We mention here that studying the improvement versions (or stability estimates) of inequalities in analysis and geometric recently have attracted lots of attentions by many mathematicians and becomes an iteresting field in mathematical research, for examples see \cite{BB10, BB11,CFMP09,FiMP09,FMP10,FMP13,FMP08} and references therein. For our purpose, it is convenient to introduce the $L_p-$mixed volume deficit, $L_p-$Firey-Brunn-Minkowski deficit and the relative side factor of $K, L\in \mK_0^n$, respectively, by
$$\de_p(K,L) = \frac{V_p(K,L)}{V(K)^{1-\frac pn} V(L)^{\frac pn}} -1,\quad \be_p(K,L) = \frac{V(K+_p L)^{\frac pn}}{V(K)^{\frac pn} + V(L)^{\frac pn}}-1,$$
and
$$\si(K,L) = \max\lt\{\frac{V(K)}{V(L)},\frac{V(L)}{V(K)}\rt\}.$$

The first main result of this paper is a bound from below of $\de_p(K,L)$ in terms of the relative asymmetry index of $K$ and $L$, defined as
$$A(K,L) = \frac{V(K\De \lam L)}{V(K)}\quad\text{with}\quad \lam = \lt(\frac{V(K)}{V(L)}\rt)^{\frac 1n},$$
where $\De$ denotes the symmetric difference of two subsets of $\R^n$. More precisely, we have the following theorem.
\begin{theorem}\label{eq:1stTheo}
Let $p> 1$ be a real number. Then for any convex bodies $K,L\in \mK_0^n$, it holds
\begin{equation}\label{eq:stabresult}
\de_p(K,L) \geq \frac{p-1}{128n^2} A(K,L)^2.
\end{equation}
\end{theorem}
From Theorem \ref{eq:1stTheo}, we derive the following improvement of \eqref{eq:extensionBM}.
\begin{theorem}\label{stabpBM}
Let $p > 1$ be a real number. Then for any convex bodies $K,L\in\mK_0^n$, it holds
\begin{equation}\label{eq:stabpBM}
\be_p(K,L) \geq \frac{p-1}{512 n^2\,\si(K,L)^{\frac pn}}\, A(K,L)^2.
\end{equation}
\end{theorem}

Theorem \ref{eq:1stTheo} and Theorem \ref{stabpBM} do not give any infomation about the below bounds of $\de_p(K,L)$ and $\beta_p(K,L)$ when $p=1$. In fact, finding the stablity version for the $L_1-$mixed volume inequality (in other word, the anisotropic isoperimetric inequality) and Brunn-Minkowski inequality are extreme difficult problems which are recently proved in \cite{FMP10,FMP08} and in \cite{FJ15}, respectively. It is worth pointing out that the order of $A(K,L)$ in Theorem \ref{eq:1stTheo} and Theorem \ref{stabpBM} is sharp. An example is given at the end of section \S3 to show this optimality.


We conclude this section by introducing a stability version for Jensen inequality. Let $(\Om,\mu)$ be an arbitrary probability measure space and let $p$ be a positive number, $p\not=1$. Jensen inequality asserts that for any nonnegative function $f\in L^p(\mu)$, it holds
\begin{equation}\label{eq:Jensen}
\frac1{p-1}\lt(\int_\Om f^p \, d\mu - \lt(\int_\Om f\, d\mu\rt)^p\rt) \geq 0.
\end{equation}
The next theorem shows that we can strengthen the inequality \eqref{eq:Jensen} by adding a remainder term which measures the deviation between $\frac f{\int_\Om f\, d\mu}$ and $1$. More precisely, we have the following theorem.
\begin{theorem}\label{stabHol}
Given $p > 0$, let us denote 
\begin{equation*}
c_p = 
\begin{cases}
\frac12 &\mbox{if $p\geq 1$}\\
\frac{(p+1)^{p+1}}{8p^{p-1}} &\mbox{if $0<p<1$.}
\end{cases}
\end{equation*}
Then for any nonnegative function $f\in L^p(\mu)\cap L^1(\mu)$, we have 
\begin{equation}\label{eq:stabHol}
\frac1{p-1}\lt(\frac{\int_\Om f^p d\mu}{\lt(\int_\Om f d\mu\rt)^p} -1\rt) \geq  c_p\lt(\int_\Om \lt|\frac f{\int_\Om f d\mu} -1\rt| d\mu\rt)^2.
\end{equation}
When $p=1$, the left hand side of \eqref{eq:stabHol} should be understood as $\int_\Om f \ln (f/ \int_\Om f d\mu)\, d\mu$ which is Shannon entropy of $f$.
\end{theorem}
Theorem \ref{stabHol} is an important ingredient in the proof of Theorem \ref{eq:1stTheo} as shown below. The interest of this theorem is that it recovers, in the case $p=1$, the famous Pinsker-Csisz\'ar-Kullback inequality for Shannon entropy which has many applications in Information theory (see \cite{Bo07,Csis67,Kull68,Pin64}). For $0< p\not=1$ and any nonnegative function $f\in L^p(\mu)$, such that $\int_\Om f\, d\mu =1$, let us introduce the $p-$Tsallis entropy of $f$ by
$$S_p(f) = \frac1{p-1}\lt(\int_\Om f^p\, d\mu -1\rt).$$
Note that $\lim\limits_{p\to 1} S_p(f) = \int_\Om f\ln f\, d\mu$ which is Shannon entropy of $f$. Theorem \ref{stabHol} hence can be seen as a generalization of Pinsker-Csisz\'ar-Kullback inequality for Tsallis entropy. We should emphasize here that there are many improvements of Jensen inequality in literature (e.g., see \cite{Dra90,Dra92,Dra94,Dra95,Dra10}), but they are very different with the one in Theorem \ref{stabHol}. Theorem \ref{stabHol} seems to be new, and maybe is of independent interest.


The rest of this paper is organized as follows. In section \S2, we give the proof of the improvement of Jensen inequality (Theorem \ref{stabHol}). In section \S3, we use Theorem \ref{stabHol} to prove Theorem \ref{eq:1stTheo} and then derive Theorem \ref{stabpBM}.

\section{Proof of Theorem \ref{stabHol}}
In this section, we give the proof of Theorem \ref{stabHol}. The proof of Pinsker-Csisz\'ar-Kullback inequality inequality (case $p=1$) can be found in \cite{Bo07}. Hence, from now on, we only consider the case $p > 0$ and $p\not=1$.  

By the density argument and the homogeneity, we can assume, without loss of generality, that $f > 0$ on $\Om$ and $\int_\Om f d\mu =1$. Let us denote $A =\{x\,:\, 0< f(x) \leq 1\}$, $t = \mu(A)$ and $a = \int_A f d\mu$, then  $a \leq t \leq 1$. 

Let us first prove that $a > 0$. Indeed, if $a=0$, then we must have $\mu(A)=0$, or equivalently $\mu(A^c)=1$. Since $f > 1$ on $A^c$, hence it holds
$$1 = \int_\Om f d\mu = \int_{A^c} f d\mu  > \mu(A^c) =1.$$
This contradiction shows that $a > 0$.

If $t= 1$ then $a=1$. This implies that $f =1$ almost everywhere on $\Om$ hence the inequality \eqref{eq:stabHol} is trivial in this case.

It remains to consider the case $a,t\in (0,1)$. If $p> 1$, it follows from H\"older inequality that
$$a^p =\lt(\int_A f d\mu\rt)^p \leq \mu(A)^{p-1}\int_A f^p d\mu,$$
or equivalently,
\begin{equation}\label{eq:a1}
\int_A f^p d\mu \geq t^{1-p} a^p.
\end{equation}
The similar argument shows that
\begin{equation}\label{eq:a2}
\int_{A^c} f^p d\mu \geq (1-t)^{1-p}(1-a)^p.
\end{equation}
If $p\in (0,1)$, the inequalities \eqref{eq:a1} and \eqref{eq:a2} change the direction by reverse H\"older inequality. Consequently, we get
$$\frac1{p-1}\lt(\frac{\int_\Om f^p d\mu}{\lt(\int_\Om f d\mu\rt)^p} -1\rt) \geq \frac{t^{1-p} a^p +(1-t)^{1-p}(1-a)^p-1}{p-1}.$$
In the other hand,
$$\int |f -1| d\mu = 2\int_A (1-f)d\mu = 2(t-a).$$
Hence, it is enough to prove that 
\begin{equation}\label{eq:enough}
\frac{t^{1-p} a^p +(1-t)^{1-p}(1-a)^p-1}{p-1}-4c_p(t-a)^2 \geq 0,\quad \forall\, t\in [a,1].
\end{equation}
Let $\psi_a(t)$ denote the left hand side of \eqref{eq:enough}. By a direct computation, we have
$$\psi_a'(t) =(1-t)^{-p}(1-a)^p -t^{-p}a^p -8c_p(t-a).$$
An application of the fundamental theorem of calculus yields
$$ (1-t)^{-p}(1-a)^p -t^{-p}a^p=p \int_{\frac at}^{\frac{1-a}{1-t}} s^{p-1} ds = p\frac{t-a}{t(1-t)} \int_0^1\lt((1-s) \frac at + s \frac{1-a}{1-t}\rt)^{p-1} ds.$$
We devide the proof of \eqref{eq:enough} into two cases according to $p> 1$ or $p\in (0,1)$.
\begin{description}
\item (i) We first consider the case $p\in (0,1)$. Since $a/t \leq (1-a)/(1-t)$ then
$$(1-t)^{-p}(1-a)^p -t^{-p}a^p \geq p\frac{(1-a)^{p-1}(t-a)}{t(1-t)^p} \geq p\frac{t-a}{t(1-t)^p}.$$
The function $t(1-t)^p$ attains its maximum at $t = 1/(p+1)$ on $[0,1]$, hence
$$(1-t)^{-p}(1-a)^p -t^{-p}a^p \geq \frac{(p+1)^{p+1}}{p^{p-1}}(t-a) = 8c_p (t-a).$$
Consequently, we have 
$$\psi_a'(t) \geq 0,\quad\forall\, t\geq a,$$
and then $\psi_a(t) \geq \psi_a(a) =0$ for any $t\geq a$.

\item (ii) Let us consider the case $p> 1$. Since $1-a \geq 1-t > 0$ then
$$(1-t)^{-p}(1-a)^p -t^{-p}a^p \geq \frac{(t-a)}{t (1-t)} \int_0^1 ps^{p-1} ds = \frac{(t-a)}{t (1-t)}\geq 4(t-a),$$
for any $t\geq a$. Consequently, we have 
$$\psi_a'(t) \geq 0,\quad\forall\, t\geq a,$$
and then $\psi_a(t) \geq \psi_a(a) =0$ for any $t\geq a$.
\end{description}
We conclude that \eqref{eq:enough} holds. The proof of Theorem \ref{stabHol} is complete.

An immediate consequence of Theorem \ref{stabHol} is an improvement of Jensen inequality for concave funtion $\phi(t) = \ln t$ as follows,
$$\ln\lt(\int_\Om f d\mu \rt) - \int_\Om \ln f\, d\mu \geq \frac18 \lt(\int_\Om \lt|\frac f{\int_\Om f d\mu}-1\rt|d\mu\rt)^2.$$
Indeed, This inequality follows by dividing both sides of \eqref{eq:stabHol} by $p$, and then letting $p$ tend to $0$.

\section{Proof of Theorem \ref{eq:1stTheo} and Theorem \ref{stabpBM}}
We start this section by giving a proof of Theorem \ref{eq:1stTheo}. The proof is an application of Theorem \ref{stabHol} to the integral representation \eqref{eq:intrepre} for the $L_p-$mixed volume. Let us go into details.

\emph{Proof of Theorem \ref{eq:1stTheo}:} By homogeneity, we can assume, without loss of generality, that $V(K)=V(L)=1$. Under these assumptions, the $L_p-$mixed volume deficit $\de_p(K,L)$ has form
$$\de_p(K,L) = V_p(K,L) -1 = \frac 1n \int_{S^{n-1}} \lt(\frac{h(L,u)}{h(K,u)}\rt)^p h_K(u) dS(K,u)-1.$$
By the integral representation \eqref{eq:intrepre}, we have
$$\frac1n \int h( L,u) dS(K,u) = V_1(K, L).$$
The $L_1-$mixed volume inequality \eqref{eq:Lpmixed} implies $V_1(K,L) \geq 1$. Let us denote $\ga = \frac1{V_1(K,L)} \leq 1$ for convenient. Theorem \ref{stabHol} leads to
\begin{align}\label{eq:aaa}
\de_p(K,L)& \geq \frac{p-1}2 V_1(K, L)^p \lt(\frac1n\int_{S^{n-1}} \lt|\frac{h(\ga L,u)}{h(K,u)} -1\rt| h(K,u) dS(K,u)\rt)^2 + V_1(K, L)^p -1\notag\\
&\geq \frac{p-1}2 \lt(\frac1n\int_{S^{n-1}} \lt|h(\ga L,u)-h(K,u)\rt| dS(K,u)\rt)^2 + p(V_1(K, L) -1).
\end{align}
In particuliar, \eqref{eq:aaa} yields
$$V_1(K, L) \leq 1 + \frac{\de_p(K,L)}p,$$
or equivalently, $\ga \geq \frac{p}{p+\de_p(K,L)}$. Denote by $K_1$ the convex hull of $K\cup \ga L$ and $K_2 = K \cap \ga L$. We then easily check that
$$h(K_1,\cdot) = \max\{h(K,\cdot), h(\ga L,\cdot)\},\quad h(K_2,\cdot) = \min\{h(K,\cdot), h(\ga L,\cdot)\}.$$
From the intergral representation \eqref{eq:intrepre}, we have
\begin{align}\label{eq:111}
\frac1n\int_{S^{n-1}} &\lt|h(\ga L,u)-h(K,u)\rt| dS(K,u)\notag\\
&=\frac12\lt(\frac1n\int_{S^{n-1}} h(K_1,u)dS(K,u)-\frac1n\int_{S^{n-1}} h(K_2,u)dS(K,u)\rt)\notag\\
&= \frac{V_1(K,K_1) -V(K,K_2)}2\notag\\
&\geq \frac{V_1(K,K_1) -V(K)}2,
\end{align}
where the inequality comes from the fact $V_1(K,K_2) \leq V_1(K,K) = V(K)$ since $K_2\subset K$. The $L_1-$mixed volume inequality implies 
\begin{equation}\label{eq:bbb}
V_1(K,K_1) -V(K) \geq V(K_1)^{\frac1n} -V(K)^{\frac1n} \geq \frac{V(K_1)-V(K)}{nV(K_1)^{\frac{n-1}n}}.
\end{equation}
Combining \eqref{eq:aaa}, \eqref{eq:111} and \eqref{eq:bbb} shows that
$$\de_p(K,L) \geq \frac{p-1}8 (V(K_1)^{\frac1n} -1)^2,$$
or equivalently,
\begin{equation}\label{eq:ccc}
V(K_1) \leq \lt(1 + 2\sqrt{\frac{2\de_p(K,L)}{p-1}}\rt)^n=:\ga_p(K,L)^{n}.
\end{equation}
Combining \eqref{eq:bbb} and \eqref{eq:ccc} leads to
\begin{equation}\label{eq:ddd}
V_1(K,K_1) -V(K) \geq \frac{V(K_1) -V(K)}{n \ga_p(K,L)^{n-1}}.
\end{equation}
Since $K\cup (\ga L\setminus K) = K \cup (\ga L) \subset K_1$, hence
\begin{equation}\label{eq:eee}
V(K_1) -V(K) \geq V(\ga L\setminus K) = V(L\setminus K) - V((L\setminus \ga L)\setminus K).
\end{equation}
Plugging \eqref{eq:eee}, \eqref{eq:ddd}, \eqref{eq:111} and the fact $\ga\leq 1$ into \eqref{eq:aaa}, we obtain
\begin{align*}
\de_p(K,L) &\geq \frac{p-1}2\lt(\frac{V(L\setminus K) - V((L\setminus \ga L)\setminus K)}{n\ga_p(K,L)^{n-1}}\rt)^2 + p(1-\ga).
\end{align*}
Using the simple inequalities $V(L\setminus K)\leq 1$ and
$$V((L\setminus \ga L)\setminus K)\leq V(L\setminus \ga L) = 1-\ga^n,$$
we can readily prove that
$$\de_p(K,L) \geq \frac{p-1}{2n^2 \ga_p(K,L)^{2(n-1)}} V(L\setminus K)^2 - \frac{p-1}{n\ga_p(K,L)^{n-1}} (1-\ga^n) + p(1-\ga).$$
Since $\ga_p(K,L) \geq 1$ and $1-\ga^n \leq n(1-\ga)$, hence we have
\begin{equation}\label{eq:fff}
\de_p(K,L) \geq \frac{p-1}{2n^2 \ga_p(K,L)^{2(n-1)}} V(L\setminus K)^2 =\frac{p-1}{8n^2 \ga_p(K,L)^{2(n-1)}} V(K\De L)^2.
\end{equation}

If $n=1$, then \eqref{eq:stabresult} is a trivial consequence of \eqref{eq:fff}. Let us consider the case $n\geq 2$. It is clear that $V(K\De L) \leq 2$. If $\de_p(K,L) \geq \frac{p-1}{32(n-1)^2}$ then 
$$\de_p(K,L) \geq \frac{p-1}{128(n-1)^2}V(K\De L)^2\geq \frac{p-1}{128n^2}V(K\De L)^2.$$
If $\de_p(K,L) \leq \frac{p-1}{32(n-1)^2}$ then $\ga_p(K,L) \leq 1 + \frac{1}{2(n-1)}$ by \eqref{eq:ccc}, hence \eqref{eq:fff} yields
$$\de_p(K,L) \geq \frac{p-1}{8n^2 \lt(1 +\frac{1}{2(n-1)}\rt)^{2(n-1)}}V(K\De L)^2\geq\frac{p-1}{24n^2}V(K\De L)^2 \geq \frac{p-1}{128n^2}V(K\De L)^2.$$
Theorem \ref{eq:1stTheo} therefore is completely proved.

We next show how to derive Theorem \ref{stabpBM} from Theorem \ref{eq:1stTheo}. The idea goes back \cite{FMP10} where the authors obtain a stability estimate for Brunn-Minkowski inequality on convex bodies from the quantitative anisotropic isoperimetric inequality. We will need the following fact  
\begin{equation}\label{eq:tamgiac}
A(K,L) \leq A(K,M) + A(M,L),\quad \forall\, K, L, M \in \mK_0^n.
\end{equation}
which is a simple consequence of the inclusion
$$K\De L \subset (K\De M) \cup (M\De L).$$

\emph{Proof of Theorem \ref{stabpBM}:} Using the identity $V_p(K,K) =V(K)$ for any convex body $K$ and the linearity of $V_p(K,L)$ in $L$ with respect to the Minkowski-Firey $L_p-$combination, we have
$$V(K+_pL) = V_p(K+_pL,K+_pL) = V_p(K+_pL,K) + V_p(K+_pL,L).$$
Theorem \ref{eq:1stTheo} yields
$$V_p(K+_pL,K) \geq V(K+_pL)^{1-\frac pn} V(K)^{\frac pn} \lt(1 +\frac{p-1}{128n^2} A(K+_pL,K)^2\rt),$$
and
$$V_p(K+_pL,L) \geq V(K+_pL)^{1-\frac pn} V(L)^{\frac pn} \lt(1 +\frac{p-1}{128n^2} A(K+_pL,L)^2\rt).$$
Summing these two inequalities and then dividing both sides of the obtained inequality by $V(K)^{\frac pn} +V(L)^{\frac pn}$, we obtain
\begin{align*}
\be_p(K,L)& \geq \frac{V(K)^{\frac pn}}{V(K)^{\frac pn}+V(L)^{\frac pn}}\frac{p-1}{128n^2}A(K+_pL,K)^2)\\
&\quad\quad + \frac{V(L)^{\frac pn}}{V(K)^{\frac pn}+V(L)^{\frac pn}} \frac{p-1}{128n^2}A(K+_pL,L)^2)\\
&\geq \frac{p-1}{256n^2\si(K,L)^{\frac pn}} (A(K+_pL,K)^2+A(K+_pL,L)^2)\\
&\geq \frac{p-1}{512n^2\si(K,L)^{\frac pn}} (A(K+_pL,K)+A(K+_pL,L))^2\\
&\geq \frac{p-1}{512n^2\si(K,L)^{\frac pn}}A(K,L)^2,
\end{align*}
where the last inequality comes from \eqref{eq:tamgiac}. This finishes the proof of Theorem \ref{stabpBM}.

We conclude this paper by proving the optimality of the order of $A(K,L)$ in \eqref{eq:stabresult} and \eqref{eq:stabpBM}. Take $K = B_2^n$ the Euclidean unit ball and $L =\ep x_0+K$ with $\ep \in (0,1)$ and $x_0\in S^{n-1}$. Then, $h(K,u) =1$ and $h(L,u) = 1+ \ep \la x_0,u\ra$ for any $u\in S^{n-1}$. It is easy to verify that 
$$A(K,L)^2 \sim \ep^2,$$ 
and
$$V_p(K,L)-V(K) = \frac1n\int_{S^{n-1}} (1+\ep \la x_0,u\ra)^pdS(K,u)-V(K) \sim \ep^2,$$
for $\ep>0$ small enough. Hence the order of $A(K,L)$ in \eqref{eq:stabresult} is sharp. 

For inequality \eqref{eq:stabpBM}, we remark that the support function of $K+_p L$ is given by
$$h(K+_pL,u) = (1 + (1 + \ep \la x_0,u\ra)^p)^{\frac1p},\quad u\in S^{n-1}.$$
By Taylor's expansion, there exist positive constants $C,c$ such that 
$$c t^2 \leq (1+(1+t)^p)^{\frac 1p}-2^{\frac 1p} -2^{\frac 1p-1} t \leq Ct^2 \leq C|t|,$$
when $t\sim 0$. Consequently, for $\ep > 0$ small enough, we have
$$ 2^{\frac 1p} + 2^{\frac1p-1} \ep \la x_0,u\ra + c\ep^2|\la x_0,u\ra|^2\leq h(K+_pL,u) \leq 2^{\frac 1p} + 2^{\frac1p-1} \ep \la x_0,u\ra + C\ep^2|\la x_0,u\ra|.$$
Denote $L_1=\{x\, :\, |\la x,u\ra| \leq |\la x_0,u\ra|^2,\, \forall\, u\in S^{n-1}\}$ and $L_2$ the segment $[-x_0,x_0]\subset \R^n$, then 
$$2^{\frac1p-1}\ep x_0 +2^{\frac1p}B_2^n + c\ep^2 L_1\subset K+_pL \subset 2^{\frac1p-1}\ep x_0 +2^{\frac1p}B_2^n + C\ep^2 L_2.$$
These inclusions yield the existence of $C_1,c_1>0$ such that 
$$c_1\ep^2 \leq V(K+_pL)^{\frac pn}-2V(B_2^n)^{\frac pn} \leq C_1 \ep^2,$$
for $\ep > 0$ small enough. Hence the order of $A(K,L)$ in \eqref{eq:stabpBM} is sharp.


\begin{thebibliography}{9999}
\bibitem{BB10}
K. Ball, and K. J. B\"or\"oczky, \emph{Stability of the Pr\'ekopa-Leindler inequality\text}, Mathematika, {\bf 56} (2010) 339-356.

\bibitem{BB11}
K. Ball, and K. J. B\"or\"oczky, \emph{Stability of some versions of the Pr\'ekopa-Leindler inequality\text}, Monatshefte Math., {\bf 163} (2011) 1-14.

\bibitem{Bo07}
V. I. Bogachev, \emph{Measure theory, Volume $1$\text}, Springer-Verlag, Berlin, 2007.



\bibitem{CFMP09}
A. Cianchi, N. Fusco, F. Maggi, and A. Pratelli, \emph{The sharp Sobolev inequality in quantitative form\text}, J. Eur. Math. Soc., {\bf 11} (5) (2009) 1105-1139.

\bibitem{Csis67}
I. Csisz\'ar, \emph{Information-type measures of difference of probability distributions and indirect observations\text}, Studia Sci. Math. Hungar., {\bf 2} (1967) 299-318.

\bibitem{Dra90}
S. S. Dragomir, \emph{An improvement of Jensen's inequality\text}, Bull. Math. Soc. Sci. Math. Roumanie, {\bf 34} (1990) 291-296.

\bibitem{Dra92}
S. S. Dragomir, \emph{Some refinements of Jensen's inequality\text}, J. Math. Anal. Appl., {\bf 168} (1992) 518-522.

\bibitem{Dra94}
S. S. Dragomir, \emph{A further improvement of Jensen's inequality\text}, Tamkang J. Math., {\bf 25} (1994) 29-36.

\bibitem{Dra95}
S. S. Dragomir, \emph{A new improvement of Jensen's inequality\text}, Indian J. Pure and Appl. Math., {\bf 26} (1995) 959-968.

\bibitem{Dra10}
S. S. Dragomir, \emph{A refinement of Jensen inequality with applications for $f-$divergence measures\text}, Taiwanese J. Math., {\bf 14} (2010) 153-164.

\bibitem{FiMP09}
A. Figalli, F. Maggi, and A. Pratelli, \emph{A refined Brunn-Minkowski inequality for convex sets\text}, Ann. Inst. H. Poincar\'e Anal. Non Lin\'eaire, {\bf 26} (2009) 2511-2519.

\bibitem{FMP10}
A. Figalli, F. Maggi, and A. Pratelli, \emph{A mass transportation approach to quantitative isoperimetric inequalities\text}, Invent. Math., {\bf 182} (2010) 167-211.

\bibitem{FMP13}
A. Figalli, F. Maggi, and A. Pratelli, \emph{Sharp stability theorems for the anisotropic Sobolev and log-Sobolev inequalities on functions of bounded variations\text}, Adv. Math., {\bf 242} (2013) 80-101.

\bibitem{FJ15}
A. Figalli, and D. Jerison, \emph{Quantitative stability for the Brunn-Minkowski inequality\text}, preprint,  	arXiv:1502.06513.

\bibitem{Firey62}
W. J. Firey, \emph{$p-$means of convex bodies\text}, Math. Scand., {\bf 10} (1962) 17-24.

\bibitem{FMP08}
N. Fusco, F. Maggi, and A. Pratelli, \emph{The sharp quantitative isoperimetric inequality\text}, Ann. of Math. (2), {\bf 168} (2008) 941-980.



\bibitem{Kull68}
S. Kullback, \emph{On the convergence of discrimination information\text}, IEEE Trans. Information Theory, IT-14, (1968) 765-766.

\bibitem{Lut91}
E. Lutwak, \emph{The Brunn-Minkowski-Firey theory I: mixed volumes and the Minkowski problem\text}, J. Diff. Geom., {\bf 38} (1993) 131-150.

\bibitem{Pin64}
M. S. Pinsker, \emph{Information and information stability of random variables and processes\text}, Translated and edited by Amiel Feinstein. Holden-Day Inc., San Francisco, Calif., 1964.

\end{thebibliography}
\end{document}